\documentclass[12pt]{amsart}
\usepackage{amssymb}
\usepackage{wasysym}

\newtheorem{thm}{Theorem}[section]
\newtheorem{lem}[thm]{Lemma}

\newtheorem{prop}[thm]{Proposition}
\newtheorem{ex}[thm]{Example}

\newtheorem{con}[thm]{Conjecture}
\newtheorem*{prob*}{Open problem}

\theoremstyle{definition}

\newtheorem{defi}[thm]{Definition}

\theoremstyle{remark}

\newtheorem{rem}[thm]{Remark}
\newtheorem*{rem*}{Remark}

\newenvironment{items}{\begin{list}{$\alph{item})$}
{\labelwidth18pt \leftmargin18pt \topsep3pt \itemsep1pt \parsep0pt}}
{\end{list}}
\newcommand{\bulit}{\item[$\bullet$]}

\newcommand{\kringel}{\mathbin{\raise1pt\hbox{$\scriptstyle\circ$}}}
\newcommand{\pkt}{\mathbin{\raise0pt\hbox{$\scriptstyle\bullet$}}}

\newcommand{\F}{\mathbb{F}}
\newcommand{\N}{\mathbb{N}}

\newcommand{\Q}{\mathbb{Q}}

\newcommand{\La}{\mathfrak{a}}

\newcommand{\Lf}{\mathfrak{f}}
\newcommand{\Lg}{\mathfrak{g}}
\newcommand{\Lh}{\mathfrak{h}}

\newcommand{\Ll}{\mathfrak{l}}

\newcommand{\Lz}{\mathfrak{z}}

\newcommand{\LA}{\mathfrak{A}}

\newcommand{\CI}{\mathcal{I}}

\newcommand{\al}{\alpha}

\newcommand{\om}{\omega}

\newcommand{\Ga}{\Gamma}

\newcommand{\ra}{\rightarrow}

\renewcommand{\phi}{\varphi}

\begin{document}

\title[Faithful representations]{Computing faithful representations \\
for nilpotent Lie algebras}

\author[D. Burde]{Dietrich Burde}
\author[B. Eick]{Bettina Eick}
\author[W. de Graaf]{Willem de Graaf}
\address{Fakult\"at f\"ur Mathematik\\
Universit\"at Wien\\
  Nordbergstr. 15\\
  1090 Wien \\
  Austria}
\email{dietrich.burde@univie.ac.at}
\address{Institut Computational Mathematics\\
TU Braunschweig\\
Pockelsstraße 14\\
31806 Braunschweig\\
Germany}
\email{beick@tu-bs.de}
\address{Dipartimento di Matematica\\
Via Sommarive 14\\
I-38050 Povo (Trento)\\
Italy}
\email{degraaf@science.unitn.it}

\date{\today}

\subjclass{Primary 17B30, 17D25}
\thanks{The third author expresses his gratitude towards the Erwin 
Schr\"odinger International Institute for Mathematical Physics for its 
hospitality in May 2008, when part of the work on the paper was carried out}

\begin{abstract}
We describe three methods to determine a faithful representation of small 
dimension for a finite-dimensional nilpotent Lie algebra over an arbitrary
field. We apply our methods in finding bounds for the smallest dimension
$\mu(\Lg)$ of a faithful $\Lg$-module for some nilpotent Lie algebras $\Lg$. 
In particular, we describe an infinite family of filiform nilpotent Lie 
algebras $\Lf_n$ of dimension $n$ over $\Q$ and conjecture that 
$\mu(\Lf_n) > n+1$. Experiments with our algorithms suggest that 
$\mu(\Lf_n)$ is polynomial in $n$.
\end{abstract}

\maketitle

\section{Introduction}

The Ado-Iwasawa theorem asserts that every finite-dimensional Lie algebra 
over an arbitrary field has a faithful finite-dimensional 
representation. A constructive proof for this theorem in characteristic 0
has been given in \cite{GRA4}. It has been implemented as an 
algorithm in the computer algebra systems GAP and Magma.

Here we consider a variation on this theme: we introduce three algorithms
for computing a faithful finite-dimensional representation for a finite-dimensional 
nilpotent Lie algebra over an arbitrary field. They take as
input a nilpotent Lie algebra $\Lg$ given by a structure
constants table and can be briefly summarized as follows.

\begin{items}
\bulit
Our first algorithm uses an action of $\Lg$ on a quotient of the universal 
enveloping algebra of $\Lg$. It is based on ideas in \cite{BU7}. 
\bulit
The second algorithm constructs a finite-dimensional faithful submodule of 
the dual of the universal enveloping algebra of $\Lg$. The resulting module 
is minimal in the sense that it has no faithful quotients or submodules. 
\bulit
The third algorithm uses a randomised method to try to construct a faithful
representation in dimension $\dim(\Lg)+1$. It uses induction on a central 
series of $\Lg$, in each step extending the representation by means of a 
cohomological construction. 
\end{items}

The methods described here are practical. In particular, all three methods 
are usually more efficient than the general algorithm given in \cite{GRA4}. 
We include a report on various applications of our GAP implementation of 
our algorithms below.

A central aim in the design of our algorithms was the idea to try to construct
a faithful representation of possibly small dimension of a given nilpotent
Lie algebra. In particular, our three algorithms usually determine 
representations which are a lot smaller than the representation obtained
by the algorithm from \cite{GRA4}. Our algorithms can therefore be used to 
determine good upper bounds for the minimal possible dimension $\mu(\Lg)$ 
of a faithful module of the nilpotent Lie algebra $\Lg$. 

The interest in the invariant $\mu(\Lg)$ is motivated, among other things, 
by problems from geometry and topology. For example, Milnor and Auslander 
studied generalizations of crystallographic groups; these are related to 
$\mu(\Lg)$ as follows: let $\Gamma$ be a finitely-generated, torsion-free 
nilpotent group of rank $n$ with real Malcev completion $G_{\Gamma}$. Then, 
if $\Gamma$ is the fundamental group of a compact complete affinely-flat
manifold, it follows that $\mu(\Lg_{\Ga})\le n+1$ for the Lie algebra 
$\Lg_{\Ga}$ of $G_{\Ga}$. 
Another motivation for studying $\mu(\Lg)$ is based on the result
that the Lie algebra $\Lg$ of a Lie group $G$ admitting a left-invariant 
affine structure satisfies $\mu(\Lg)\le n+1$. 
In the context of these two results Milnor conjectured that any solvable 
Lie group should admit a left-invariant affine structure. The algebraic 
formulation of this conjecture implies that $\mu(\Lg) \le \dim (\Lg)+1$ should
hold for the Lie algebra $\Lg$ of $G$.

There are counter-examples known to Milnor's conjecture. Indeed, 
there are infinitely many filiform nilpotent Lie algebras of dimension 
$10$ which do not have any faithful module of dimension $n+1$. We refer to 
\cite{BU13} for details and background.

Here we describe an infinite family of Lie algebras $\Lf_n$ of dimension
$n$ for $n \geq 13$ and we use our algorithms to study the invariant
$\mu(\Lf_n)$ for these Lie algebras. We conjecture that these 
Lie algebras do not have a faithful representation of dimension $n+1$. But
our experiments suggest that $\mu(\Lf_n)$ is polynomial in $n$ for these
Lie algebras.


\section{Using quotients of the universal enveloping algebra}
\label{sect:quot}

Let $\Lg$ be a finite-dimensional nilpotent Lie algebra over an arbitrary
field. By $\Lg^m$, $m\geq 1$, we denote the terms of the lower central 
series of $\Lg$. If $x_1, \ldots, x_d$ is a basis of $\Lg$, then the formal 
products $x_1^{\alpha_1} \cdots x_d^{\alpha_d}$ with $\alpha_i \in \N$
form a basis of the universal enveloping algebra $U(\Lg)$. 

We define the {\em weight} $wgt(x)$ of an element $x \in \Lg$ as the 
maximal $m$ with $x \in \Lg^{m}$. The {\em weight} of a basis element 
of $U(\Lg)$ is then defined by 
\[ wgt(x_1^{\alpha_1} \cdots x_d^{\alpha_d}) 
                      = \sum_{i=1}^d \alpha_i \, wgt(x_i) \]

Let $U^m(\Lg) = \langle x^\alpha \mid wgt(x^\alpha) \geq m \rangle$ the
ideal in $U(\Lg)$ generated by all basis elements of weight at least $m$ for
some $m \geq 1$. The following theorem is proved in \cite{BIR},\cite{BU7}.

\begin{thm}\label{thm:1}
If $\Lg$ is nilpotent of class $c$, then $\Lg$ acts faithfully on 
$U(\Lg)/U^{c+1}(\Lg)$ by multiplication from the left. If the considered
basis of $\Lg$ contains bases for $\Lg^m$ for every $m \geq 1$, then the
resulting representation has the dimension
\[ \nu(d, c) = \sum_{j=0}^c \binom{d-j}{c-j} p(j), \]
where $p(j)$ is the number of partitions of $j$ with $p(0)=1$.
\end{thm}

This theorem yields a straightforward algorithm to construct a faithful 
module for $\Lg$. We consider a basis of $\Lg$ which contains bases of
$\Lg^m$ for every $m$. We form the space $V$ spanned by all basis elements 
of weight at most $c$ in $U(\Lg)$. An element $x \in \Lg$ acts on $V$ by left
multiplication, where we treat any element of weight at least $c+1$ as zero.
We demonstrate this algorithm in the following example.

\begin{ex}\label{exa1}
Let $\Lg$ be the $3$-dimensional Heisenberg Lie algebra spanned by $x,y,z$ 
with non-zero bracket $[x,y]=z$. The nilpotency class of $\Lg$ is $2$. We 
form the space spanned by the basis elements of $U(\Lg)$ of weight at most
2. These are 
$$1,x,y,z,x^2,xy,y^2.$$
Thus we obtain a 7-dimensional representation for $\Lg$. It is straightforward
to determine this representation explicitly by computing actions. For example,
$y\cdot x = yx = xy-z$ and $x\cdot xy = x^2y =0$. 
\end{ex}

The modules resulting from this simple and very efficient algorithm still have
rather large dimension. In the remainder of this section, we describe two 
methods to determine a module of smaller dimension from this given module.
As a first step, we note that a nilpotent Lie algebra $\Lg$ acts faithfully on 
a module if and only if its center $Z(\Lg)$ acts faithfully. Thus if $I$ is
an ideal in $U(\Lg)$ such that $I \cap Z(\Lg) = 0$, then $L$ acts faithfully
on $U(\Lg)/I$. 

For our first method we assume that the considered basis of $\Lg$ additionally
contains a basis for the center $Z(\Lg)$. We wish to determine an ideal $I$
in $U(\Lg)$ which has possibly small codimension and satisfies that $I \cap 
Z(\Lg) = 0$. Let $B$ be the set of all basis elements of weight at least 
$c+1$ in $U(\Lg)$ and initialise $I = \langle B \rangle$. We now iterate 
the following procedure: let $a$ be one of the finitely many basis elements 
of $U(\Lg)$ not contained in $I$ and not contained in $Z(\Lg)$. If $xa \in 
I$ for all $x \in \Lg$, then we add $a$ to $B$ and thus enlarge $I$ without
destroying the property $I \cap Z(\Lg) = 0$. This approach usually yields
a rather small dimensional faithful representation of $\Lg$. We demonstrate
this in the following example.
 
\begin{ex}\label{exa:ridmons}
We continue Example \ref{exa1}. We initialise $B$ as the basis elements of
weight at least $c+1$. Note that $Z(\Lg) = \langle z \rangle$. Thus we 
consider for $a$ the elements $x, y, x^2, xy, y^2$. Of those, the elements
$x^2$, $xy$ and $y^2$ satisfy the condition of the algorithm and thus
we move these into $B$. Now also $y$ satisfies the condition and we also
move $y$ into $B$. Now $I$ is an ideal of codimension $3$ in $U(\Lg)$ and
hence we obtain a 3-dimensional faithful representation of $\Lg$.
\end{ex}

So we now have a procedure to construct a possibly small set of monomials,
that yields a faithful $\Lg$-module: we first carry out the procedure of 
Theorem \ref{thm:1}, and then the procedure outlined above, to make the 
set of monomials smaller. In the sequel this algorithm will be called 
{\em Regular}.

A second method, that can be tried for any faithful $\Lg$-module $V$, is to 
perform the following algorithm on $V$.  

\begin{enumerate}
\item Compute the space $S = \{ v\in V\mid x\cdot v=0 \text{ for all }
x\in \Lg\}$.
\item Compute the space $C = \{ x\cdot v\mid v\in V, x\in Z(\Lg)\}$.
\item Set $M = S\cap C$ and let $W$ be a complement to $M$ in $S$.
\item If $W=0$ then the algorithm stops, and the output is $V$. Otherwise,
set $V := V/W$, and return to $(1)$.
\end{enumerate}

We note that any subspace of $S$ is a $\Lg$-submodule of $V$. Therefore
the quotient $V/W$ is a $\Lg$-module. Let $U$ be a complement to $W$ in $V$
such that $C\subset U$. Then $x\cdot V\subset U$ for all $x\in Z(\Lg)$. So
since $Z(\Lg)$ acts faithfully on $V$ it acts faithfully on $V/W$. Hence
 $V/W$ is a faithful $\Lg$-module.

\begin{ex}
We consider the module of Example \ref{exa1}. Here we get
\begin{align*}
S&= \langle z,x^2,xy,y^2\rangle\\
C&= \langle z\rangle \\
W&= \langle x^2, xy, y^2\rangle.
\end{align*}
After taking the quotient we get a module spanned by (the images of)
$1,x,y,z$. For this module we can perform the algorithm again. We get
$S=\langle y,z\rangle$, $C=\langle z\rangle$, $W=\langle y\rangle$. 
So we end up with a faithful module of dimension $3$.
\end{ex}

Now the complete algorithm to construct a small-dimensional faithful
$\Lg$-module consists of first performing the algorithm Regular, 
followed by the quotient procedure described above. This algorithm will
be called {\em Quotient}.


\section{Using the dual of the universal enveloping algebra}\label{sect:dual}

Now $\Lg$ acts on the dual
$U(\Lg)^*$ by $x\cdot f(a) = f(-xa)$. Let $z_1,\ldots,z_r$ be a basis 
of the center of $\Lg$, which we assume to be a subset of the basis
$x_1,\ldots,x_n$. Let $\psi_i \in U(\Lg)^*$ for $1\leq i\leq r$ be defined by
$\psi_i(z_i) = 1$ and $\psi_i(a) = 0$ for any PBW-monomial not equal to $z_i$
(note that this definition depends on the choice of basis of $\Lg$).

Let $\bar{\phantom{a}}: U(\Lg)\to U(\Lg)$ be the antiautomorphism induced
by $\bar{x}=-x$ for $x\in \Lg$. Then for $a,b\in U(\Lg)$, $f\in U(\Lg)^*$ we
get $a\cdot f(b) = f(\bar{a}b)$. (In other words, $\bar{\phantom{a}}$ is
the antipode of $U(\Lg)$.) Note that $\bar{\bar{a}}=a$.

\begin{thm}\label{thm:dual}
Let $V$ be the $\Lg$-submodule of $U(\Lg)^*$ generated by $\psi_1,\ldots,\psi_r$.
Then $V$ is a faithful finite-dimensional $\Lg$-module. Moreover, $V$ has
no faithful $\Lg$-submodules, nor has it faithful quotients.
\end{thm}

\begin{proof}
For $k\geq 1$ let $U^k(\Lg)$ be as in Theorem \ref{thm:1}.
Let $W = \{ f \in U(\Lg)^* \mid f(U^{c+1}(\Lg))=0\}$. Then $W$ is finite-dimensional
(since $U^{c+1}(\Lg)$ has finite codimension), and a $\Lg$-submodule of 
$U(\Lg)^*$ (since $U^{c+1}(\Lg)$ is an ideal). Now $V\subset W$; hence $V$
is finite-dimensional. Let $z = \sum_i \mu_i z_i$ be an element of the
center of $\Lg$. Then $z\cdot \psi_i(1) = -\psi_i(z) =-\mu_i$. So $z$ acts as
zero if and only if all $\mu_i$ are zero. So $V$ is a faithful module.\par
Let $\psi_0 \in U(\Lg)^*$ be defined by $\psi_0(1)=1$ and $\psi_0(a)=0$ for
PBW-monomials $a\neq 1$. Then $\psi_0 = -z_1\cdot \psi_1$, so $\psi_0\in V$.
Let $f\in V$, and suppose there is a PBW-monomial $a\neq 1$ such that
$f(a)\neq 0$. Then $\bar{a}\cdot f(1) = f(a) \neq 0$. The conclusion is
that $\psi_0$ spans the space of elements that are killed by $\Lg$. 
Set $M= V/\langle\psi_0\rangle$. Then $M$ is not a faithful $\Lg$-module.
Indeed, $V$ has a basis consisting of $a\cdot \psi_i$ for $1\leq i\leq r$,
and various PBW-monomials $a$. Let $z$ lie in the center of $\Lg$; then
for all PBW-monomials $b$ we get $z\cdot (a\cdot \psi_i) (b) = 
\psi_i(\bar{a}b\bar{z})$ which is zero unless $a=b=1$ (note that 
$\bar{z}=-z$ also lies in the center of $\Lg$). Furthermore, $z\cdot \psi_i
= -\mu_i \psi_0$ (where $z= \sum_j \mu_j z_j$). It follows that the center
of $\Lg$ acts trivially on $V/\langle \psi_0\rangle$. In particular it is
not a faithful $\Lg$-module. Now, since every $\Lg$-submodule of $V$ must
contain $\psi_0$, it follows that $V$ has no faithful quotients.\par
Let $M\subset V$ be a faithful $\Lg$-submodule. Let $b_{ij}\in U(\Lg)$ be 
such that $\{ b_{ij}\psi_i \}$ is a basis of $V$. We assume that $b_{i1}=1$,
and that the $b_{ij} \in \Lg U(\Lg)$ if $j>1$ (i.e., they have no constant term).
Then $z_i\cdot b_{i1}\psi_i = -\psi_0$ and $z_i\cdot b_{kj}\psi_{k}=0$ if
$k\neq i$ or $j>1$. So since the center acts faithfully on $M$ it follows
that $M$ contains elements of the form
\[ 
\varphi_i=\psi_i + \sum_{\substack{j>1 \\ 1\leq k\leq r}} c_{kj}b_{kj}\psi_k,
\]
for $1\leq i\leq r$. (Here $c_{kj}$ are coefficients in the ground field.)
Now we introduce a weight function on $U(\Lg)^*$. 
For $k\geq 0$ set $F_k = \{f\in U(\Lg)^* \mid f(U^k(\Lg)) = 0\}$, where 
$U^k(\Lg)$ is as in Theorem \ref{thm:1}. Then $0=F_0\subset F_1\subset \cdots $.
We set $wgt(f) = k$ if $f\in F_k$ but $f\not\in F_{k-1}$. (For example: 
$wgt(\psi_0)=1$.) Let $f\in  U(\Lg)^*$ have weight $k$, and let $a\in U(\Lg)$ 
with $wgt(a)=t$; then a calculation shows that $wgt(a\cdot f) \leq k-t$.
Hence $b_{ij}\varphi_i$ is equal to $b_{ij}\psi_i$ plus a sum of functions of
smaller weight. So if we order the $b_{ij}\psi_i$ according to weight, and
express the $b_{ij}\varphi_i$ on the basis $b_{ij}\psi_i$ we get a triangular
system. We conclude that the $b_{ij}\varphi_i$ are linearly independent. 
Hence $\dim M = \dim V$ and $M=V$. 
\end{proof}

The algorithm based on this theorem is straightforward. We illustrate it
with an example.

\begin{ex}
Let $\Lg$ be the Lie algebra of Example \ref{exa1}. For a monomial
$a\in U(\Lg)$ we denote by $\psi_a$ the element of $U(\Lg)^*$ that takes the
value $1$ on $a$, and zero on all other monomials. We compute a basis of
the submodule of $U(\Lg)^*$ generated by $\psi_z$. We have $x\cdot \psi_z(a) =
-\psi_z(xa) =0$ for all monomials $a$. Secondly, $y\cdot\psi_z(x) = 
\psi_z(-yx) = \psi_z(-xy+z) = 1$. So we get $y\cdot \psi_z = \psi_x$.
Furthermore, $z\cdot\psi_z = -\psi_1$ and $\Lg\cdot\psi_1=0$, $x\cdot \psi_x = 
-\psi_1$, $y\cdot\psi_x = z\cdot \psi_x =0$. So the result is a $3$-dimensional
$\Lg$-module.
\end{ex}

\begin{rem}
Since we are working in the dual of an infinite-dimensional space it is
not immediately clear how to implement this algorithm. We proceed as follows.
Let $V$ be as in Theorem \ref{thm:dual}.
{}From the proof of Theorem \ref{thm:dual} it follows that $f(U(\Lg)^{c+1})=0$
for all $f\in V$. In other words, for all monomials $a$ with $wgt(a)\geq c+1$
and all $f\in V$ we have $f(a)=0$. It follows that we can represent an $f\in V$
by the vector containing the values $f(a)$, where $a$ runs through the
monomials of weight $\leq c$. This enables us to perform the operations of
linear algebra (testing linear dependence, constructing bases of subspaces
and so on) with the elements of $V$. Furthermore, we can compute the action
of elements of $\Lg$ on $V$.\par
The disadvantage of this approach is that the number of monomials that
have to be considered can be very large. So in the same way as in the
previous section we try to throw some monomials away.
Let $A$ be the set of monomials relative to which we represent the elements
of $U(\Lg)^*$. At the start this will be the set of monomials of weight
$\leq c$. Let $B$ be the set of all other monomials. So at the outset $B$ spans
a left ideal of $U(\Lg)$ and $f(b)=0$ for all $f\in V$ and $b\in B$. 
We move elements from $A$ to $B$,
without changing this last property. Let $a\in A$ be such that 
$a\not\in Z(\Lg)$ and $xa$ is a linear combination of elements
of $B$ for all $x\in \Lg$. Then we claim that $f(a)=0$ for all $f\in V$.
In order to see this we use the basis $\{ b_{ij}\psi_i\}$ used in the 
proof of Theorem \ref{thm:dual}. If $j=1$ then $b_{ij}\psi_i(a) = \psi_i(a)
=0$ as $a\not\in Z(\Lg)$. If $j>1$ then $b_{ij}\in \Lg U(\Lg)$ and hence
$\bar{b}_{ij}a$ is a linear combination of elements in $B$. Hence
$b_{ij}\psi_i(a) = 0$. Also the span of $B$ along with $a$ continues to be 
a left ideal. We conclude that we can move $a$ from $A$ to $B$.
We continue this process until we we do not find such monomials any more.
The resulting set is usually a lot smaller than the initial one.
\end{rem}

We note that the procedure described in the previous remark is exactly
the same as the second phase of the algorithm Regular (see Section
\ref{sect:quot}). So we first perform the algorithm Regular, and use
the resulting set of monomials to represent elements of the dual of
$U(\Lg)$. The resulting algorithm is called {\em Dual}.


\section{Affine representations at random}\label{sect:affine}

Let $\Lg$ be a nilpotent Lie algebra of dimension $d$. A homomorphism
$\rho : \Lg\to \La\Lf\Lf (K^d)\subseteq \Lg\Ll_{d+1}(K)$ into the Lie algebra 
of affine transformations 
\[
\La\Lf\Lf (K^d)\simeq \Lg\Ll (K^d)\ltimes K^d
\]
is called an {\it affine representation} of $\Lg$.
In this section we describe a method that tries to determine a faithful 
affine representation of $\Lg$ of dimension $d+1$. If the method succeeds, 
then it returns such a faithul representation of dimension $d+1$. 
However, it may also happen that the method fails and does not return a 
representation.
Also, it is worth noting that the algorithm uses random methods and hence
different runs of the algorithm may produce different results.

The method uses induction on a central series in $\Lg$. Thus we assume by
induction that we have given a central ideal $I$ in $\Lg$ with $\dim(I) = 1$
and a faithful affine representation
\[ \rho : \Lg/I \ra M_d(K). \]
Let $\{ a_1, \ldots, a_d \}$ be a basis of $\Lg$ with $I = \langle a_d 
\rangle$ and let $M_i = \rho(a_i+I)$ for $1 \leq i \leq d-1$. We assume
that every $M_i$ is a lower triangular matrix. Clearly, we can readily
extend $\rho$ to an affine representation of $\Lg$ with $\rho(a_i) = M_i$
for $1 \leq i \leq d$ where we set $M_d=0$ (so that $\rho(a_d) =0$). 
This extended representation has kernel $I$. 

Our aim is to extend $\rho$ to a faithful affine representation
\[ \psi : \Lg \ra M_{d+1}(K)\] 
such that 
\[ \psi(a_i) = \left( \begin{array}{cc}
    M_i & v_i \\ 0 & 0 \end{array} \right) \mbox{ for } 1 \leq i \leq d, \]
for certain vectors $v_i \in K^d$. The following lemma shows that the
possible values for $v_i$ can be determined using a cohomology computation.
Recall that 
\[
Z^1(\Lg, K^d) = \{ \nu : \Lg \ra K^d \mbox{ linear} \mid 
\nu([x,y]) = \rho(x) \nu(y) - \rho(y) \nu(x)\}
\]
is the space of 1-cocycles with values in the $\rho(\Lg)$-module $K^d$.

\begin{lem}
$\psi$ is a representation of $\Lg$ if and only if $v_i = \delta(a_i)$ for
$1 \leq i \leq d$ for some $\delta \in Z^1(\Lg, K^d)$.
\end{lem}

\begin{proof}
Let $\delta \in Z^1(\Lg, K^d)$ with $\delta(a_i) = v_i$. The linearity of
$\delta$ implies that $\psi$ is linear. The defining condition for maps
in $Z^1(\Lg, K^d)$ implies that $\psi$ is a Lie algebra representation.
The converse follows with similar arguments.
\end{proof}

Note that $Z^1(\Lg, K^d)$ is a vector space over $K$ and can be computed 
readily using linear algebra methods. The computation of $Z^1(\Lg, K^d)$
allows to describe all affine representations of $\Lg$ extending $\rho$. It
remains to determine the faithful representation among these. 

\begin{lem}
$\psi$ is faithful if and only if $v_{d+1} \neq 0$.
\end{lem}

\begin{proof}
If $\psi$ is faithful, then $v_{d+1} \neq 0$. Conversely, suppose that
$v_{d+1} \neq 0$. As $\rho$ is faithful, it follows that $ker(\psi) 
\subseteq I$. As $v_{d+1} \neq 0$, we find that $ker(\psi) = 0$.
\end{proof}

These ideas can be combined to the following algorithm.

\begin{enumerate}
\item
Choose a central series $\Lg = \Lg_0 > \Lg_1 > \ldots > \Lg_d > \Lg_{d+1} 
= 0$ of ideals in $\Lg$ such that $\dim(\Lg_i/\Lg_{i+1}) = 1$.
\item
by induction, extend a faithful affine representation from $\Lg/\Lg_i$
to $\Lg/\Lg_{i+1}$:
\begin{itemize}
\item[-] Compute $Z^1(\Lg/\Lg_{i+1}, K^i)$.
\item[-] Choose a $\delta \in Z^1(\Lg/\Lg_{i+1}, K^i)$ with $\delta(a_i)\neq 0$.
\item[-] If no such $\delta$ exists, then return fail.
\item[-] If $\delta$ exists, then extend $\rho$ to $\Lg/\Lg_{i+1}$.
\end{itemize}
\end{enumerate}

If $\Lg$ has a faithful affine representation of dimension $d+1$, then this 
algorithm can in principle find it. However, it may be that a ``wrong''
choice of a $\delta$ at a certain step may cause the algorithm to fail at
a later step.

The algorithm is based on linear algebra only and hence is very effective.
It often suceeds in finding a faithful representation in dimension $d+1$ if it 
exists.


\section{A series of filiform nilpotent Lie algebras}

Let $K$ be a field of characteristic zero. In this section we define a 
filiform Lie algeba $\Lf_n$ in each dimension $n\ge 13$ having interesting 
properties concerning Lie algebra cohomology, affine structures and faithful 
representations. In fact, we believe that the algebras $\Lf_n$ are counter 
examples to the conjecture of Milnor mentioned in the introduction, i.e., 
that $\mu (\Lf_n)\ge n+2$ holds. Hence it is interesting to compute the 
invariants $\mu(\Lf_n)$. 

Define an index set $\CI_n$ by
\begin{align*}
\CI_n^0 &=\{(k,s)\in \N \times \N \mid 2 \le k \le [n/2],\,
2k+1 \le s \le n \},\\
\CI_n& =\begin{cases}
\CI_n^0 & \text{if $n$ is odd},\\
\CI_n^0 \cup \{(\frac{n}{2},n)\} & \text{if $n$ is even}.
\end{cases}
\end{align*}

Now fix $n\ge 13$. We define a filiform Lie algebra $\Lf_n$ of dimension $n$ 
over $K$ as follows. For $(k,s)\in \CI_n$ let $\al_{k,s}$ be a set of 
parameters, subject to the following conditions: 
all $\al_{k,s}$ are zero, except for the following ones:

\begin{align*}
\al_{\ell,2\ell+1} & = \frac{3}{\binom{\ell}{2}\binom{2\ell-1}{\ell-1}}, \quad 
\ell=2,3,\ldots ,\lfloor \textstyle{\frac{n-1}{2}}\rfloor , \\[0.3cm]
\al_{3,n-4} & =1, \\[0.3cm]
\al_{4,n-2} & = \frac{1}{7}+\frac{10}{21}\frac{(n-7)(n-8)}{(n-4)(n-5)}, \\[0.3cm]
\al_{4,n} & =
\begin{cases}
\frac{22105}{15246}, & \text{ if $n=13$, } \\
0  & \text{ if $n \ge 14$,}
\end{cases}
\end{align*}
and 
\begin{align*}
\al_{5,n} & = \frac{1}{42}-\frac{70(n-8)}{11(n-2)(n-3)(n-4)(n-5)}+
\frac{25}{99}\frac{(n-6)(n-7)(n-8)}{(n-2)(n-3)(n-4)} \\[0.3cm]
          & \hspace{1.05cm} + \frac{5}{66} \frac{(n-5)(n-6)}{(n-2)(n-3)}-  
\frac{65}{1386} \frac{(n-7)(n-8)}{(n-4)(n-5)}.\\[0.3cm]
\end{align*}

Let $(e_1,\ldots ,e_n)$ be a basis of $\Lf_n$ and define the Lie brackets
as follows: 

\begin{align*}
[e_1,e_i] & =e_{i+1}, \quad i=2,\dots ,n-1 \\[0.3cm]
[e_i,e_j] & =\sum_{r=1}^n\biggl(\;\sum_{\ell=0}^{\lfloor \frac{j-i-1}{2}\rfloor} (-1)^\ell
\binom{j-i-\ell-1}{\ell}\al_{i+\ell,\, r-j+i+2\ell+1}\biggr)e_r,
 \quad 2 \le i<j \le n.
\end{align*}

In order to show that this defines a Lie bracket we need the following 
lemma which follows from the Pfaff--Saalsch\"utz sum formula:

\begin{lem}\label{Pfaff}
We have the following identities for all $n\ge 13$:\\[0.3cm]
\begin{align*}
\sum_{\ell=3}^{\lfloor \frac{n-1}{2}\rfloor}(-1)^{\ell-1}\binom{n-\ell-5}{\ell-2}
\al_{\ell,2\ell+1} & = \frac{(n-7)(n-8)}{(n-4)(n-5)}, \\[0.3cm]
\sum_{\ell=5}^{\lfloor \frac{n-1}{2}\rfloor}(-1)^{\ell}\binom{n-\ell-5}
{\ell-4}\al_{\ell,2\ell+1} & = -\frac{1}{70}+\frac{12(n-8)}{(n-2)(n-3)(n-4)(n-5)}, \\[0.3cm]
\sum_{\ell=3}^{\lfloor \frac{n-1}{2}\rfloor}(-1)^{\ell}\binom{n-\ell-3}{\ell-2}
\al_{\ell,2\ell+1} & = -\frac{(n-5)(n-6)}{(n-2)(n-3)}.\\[0.3cm]
\end{align*}
\end{lem}

\begin{prop}
The Jacobi identity is satisfied, so that $\Lf_n$ is a Lie algebra for any $n\ge 13$.
\end{prop}

\begin{proof}
Let $n\ge 14$ and choose the parameters $\al_{k,s}$ as follows.
Consider $\al_{k,2k+1}$, $k=3,\ldots , \lfloor \frac{n-1}{2}\rfloor$ and
$\al_{4,n-2}, \al_{5,n}$ as free variables. Let the remaining parameters be
zero, except for $\al_{2,5}=1$, $\al_{3,7}\neq 0$ and $\al_{3,n-4}=1$.
The Jacobi identity is equivalent to a system of polynomial equations in the
free parameters. First we obtain the equation $\al_{3,7}(10\al_{3,7}-\al_{2,5})=0$,
so that $\al_{3,7}=\frac{1}{10}$. More generally we see that
\begin{align*}
(\ell-1)\cdot \al_{\ell,2\ell+1} & = (4\ell+2)\cdot \al_{\ell+1,2\ell+3},  \quad 
\ell=2,3,\ldots ,\lfloor \textstyle{\frac{n-1}{2}}\rfloor .
\end{align*}
This implies the given explicit formula for all $\al_{\ell,2\ell+1}$.
Secondly we obtain 

\begin{align*}
\al_{4,n-2} & = \frac{\al_{4,9}}{\al_{3,7}}+  \frac{\al_{4,9}}{3\al_{3,7}^2}
\sum_{\ell=3}^{\lfloor \frac{n-1}{2}\rfloor}(-1)^{\ell-1}\binom{n-\ell-5}{\ell-2}
\al_{\ell,2\ell+1},
\end{align*}

\begin{align*}
\al_{5,n} & = \frac{1}{\al_{4,9}+\al_{3,7}-2\al_{2,5}}\left( 
-4\al_{4,9}+ \sum_{\ell=5}^{\lfloor \frac{n-1}{2}\rfloor}(-1)^{\ell}
\binom{n-\ell-5}{\ell-4}\al_{\ell,2\ell+1}\right) \\[0.5cm]
         & +  \frac{1}{\al_{4,9}+\al_{3,7}-2\al_{2,5}} \left(\al_{4,n-2}\left(
13\al_{4,9}+\sum_{\ell=3}^{\lfloor \frac{n-1}{2}\rfloor}(-1)^{\ell}\binom{n-\ell-3}{\ell-2}
\al_{\ell,2\ell+1} \right)\right). \\[0.3cm]
\end{align*}

This amounts to the given formulas in the definition of $\Lf_n$, if we substitute
the identities from Lemma $\ref{Pfaff}$. Conversely this also shows that the Jacobi 
identity is satisfied if the free parameters are given in this way. \\
For $n=13$ there is one difference. The parameter $\al_{4,n}$ coincides
with the parameter $\al_{4,13}$, which is given by
\[
\al_{4,13}=\frac{\al_{3,9}(-\al_{5,13}+6\al_{4,11}-5\al_{3,9})}{\al_{3,7}+2\al_{2,5}},
\]
and cannot be chosen to be zero.
For $n\ge 14$ the choice $\al_{4,n}=0$ is consistent with the Jacobi identity.
\end{proof}

\begin{ex}
The parameters for $\Lf_{13}$ are given as follows:
\begin{align*}
\al_{2,5} & = 1,\; \al_{3,7}=\frac{1}{10},\;\al_{4,9}= \frac{1}{70},\;
\al_{5,11}= \frac{1}{420}, \; \al_{6,13} = \frac{1}{2310},\\[0.3cm]
\al_{3,9} & = 1,\; \al_{4,11}=  \frac{43}{126},\; \al_{4,13}=\frac{22105}{15246},
\al_{5,13} = \frac{313}{3388}. \\
\end{align*}
\end{ex}

The algebras $\Lf_n$ belong to the family of filiform Lie algebras $\LA_n^2(K)$ 
defined in \cite{BU13}. Let us recall the following definition.

\begin{defi}
Let $\Lg$ be a filiform nilpotent Lie algebra
of dimension $n$. A $2$--cocycle $\om\in Z^2(\Lg,K)$
is called {\it affine}, if $\om \colon \Lg \wedge \Lg \rightarrow K$ does
not vanish on $\Lz (\Lg) \wedge \Lg$. A class $[\om]\in H^2(\Lg,K)$ is called
affine if every representative is affine.
\end{defi}

The cohomology class  $[\om]\in H^2(\Lg,K)$ of an affine $2$-cocycle $\om$ is affine 
and nonzero. If a filiform Lie algebra $\Lg$ of dimension $n\ge 6$ has second Betti number 
$b_2(\Lg)=2$, then there exists no affine cohomology class. \\
We have shown in \cite{BU13} that a filiform Lie algebra $\Lg$ which has an affine
cohomology class, admits a central extension
\begin{equation*}
0 \rightarrow \La \xrightarrow{\iota} \Lh \xrightarrow{\pi}
\Lg \rightarrow 0
\end{equation*}
with some Lie algebra $\Lh$ and $\iota(\La)=\Lz(\Lh)$, and has an affine structure.
In particular, such a Lie algebra has a faithful representation of dimension $n+1$. \\
We can conclude from the results in \cite{BU13} that the Lie algebras 
$\Lf_n$ do {\it not} have an affine structure arising this way.

\begin{prop}
The algebras $\Lf_n$, $n\ge 13$ have second Betti number $b_2(\Lf_n)=2$. Hence there
exists no affine cohomology class $[\om]\in H^2(\Lg,K)$. In particular there is no
central Lie algebra extension as above.
\end{prop}

For Lie algebras in $\LA_n^2(K)$ the second Betti number is $3$ or $2$, depending
on whether a certain polynomial identity $\al_{3,n-4}=P$ in the free parameters
does hold or does not hold.
For $\Lf_n$ we have chosen the parameters in such a way that $P\equiv 0$ and $\al_{3,n-4}=1$.
This implies that $b_2(\Lf_n)=2$. \\

It follows that a very natural way to obtain a faithful representation of dimension
$n+1$ does not work. In fact, we believe that there is no such representation at all
for these algebras:

\begin{con}\label{conject}
The Lie algebras $\Lf_n$, $n\ge 13$ do not have any faithful representation of 
dimension $n+1$, i.e., $\mu(\Lf_n)\ge n+2$.
\end{con}

For $n=13$ a very complicated analysis of possible faithful representations seems to
confirm this conjecture. In general our methods are not sufficient to prove this
for all $n\ge 14$. Even more difficult of course is the determination of $\mu(\Lf_n)$.

\section{Practical experiences}

We implemented all the algorithms described above in the computer algebra
system GAP. In this section we report on the application of these 
implementations to various examples. 
From Section \ref{sect:quot} we have the algorithms Regular and
Quotient. From Section \ref{sect:dual} we have the algorithm Dual.
Finally the algorithm of Section \ref{sect:affine} is called {\em Affine}.

In all our experiments Quotient and Dual returned faithful representations
of the same dimension (with Dual being slightly faster). This is illustrated
in Table \ref{Table3}. We believe that 
there must be an intrinsic reason for this to happen, such as one  module
being the dual of the other. But we have no proof of that. 
We only exhibit the results of Dual
in Tables \ref{Table1} and \ref{Table2}, noting that the results for
Quotient are similar in all cases.

All computations were done on a 2GHz processor with 1GB of memory for 
{\sf GAP}. 

\subsection{Upper triangular matrix Lie algebras}

The upper triangular matrices in $M_n(\F)$ form a nilpotent Lie algebra
$U_n(\F)$ with $n-1$ generators and class $n-1$. We applied our algorithms 
to some Lie algebras of this type. The results are recorded in Table 1.

\vspace*{0.2cm}
\begin{table}[htb]
\begin{center}
\begin{tabular}{|c|c|c|c|c|c|c|c|c|}
\hline
$n$ & $\F$ & $\dim(U_n(\F))$ 
   & \multicolumn{2}{c|}{Regular} 
   & \multicolumn{2}{c|}{Dual} 
   & \multicolumn{2}{c|}{Affine}\\
\hline
& & & time & $\dim$ & time & $\dim$ & time & $\dim$\\
\hline
4 & $\F_2$ & 6  & 0.0 & 7 & 0.1 & 5 & 0.0 & 7 \\
5 & $\F_2$ & 10 & 0.25 & 15 & 0.3 & 11 & 0.3 & 11 \\
6 & $\F_2$ & 15 & 3.4 & 35 & 3.6 & 17 & 3.5 & 16 \\
7 & $\F_2$ & 21 & 65 & 79 & 66 & 35 & 45 & 22 \\
\hline
4 & $\F_3$ & 6  & 0.0 & 7 & 0.0 & 5 & 0.0 & 7 \\
5 & $\F_3$ & 10 & 0.2 & 15 & 0.3 & 11 & 0.3 & 11 \\
6 & $\F_3$ & 15 & 3.4 & 35 & 3.6 & 17 & 3.7 & 16 \\
7 & $\F_3$ & 21 & 65 & 79 & 67 & 35 & 46 & 22 \\
\hline
4 & $\Q$ & 6  & 0.0 & 7 & 0.0 & 5 & 0.0 & 7 \\
5 & $\Q$ & 10 & 0.2 & 15 & 0.3 & 11 & 0.3 & 11 \\
6 & $\Q$ & 15 & 3.0 & 35 & 3.2 & 17 & 3.6 & 16 \\
7 & $\Q$ & 21 & 66 & 79 & 67 & 35 & 45 & 22 \\
\hline
\end{tabular}
\vspace*{0.2cm}
\caption{Running times (in seconds) for $U_n(\F)$.}\label{Table1}
\end{center}
\end{table}

Table 1 exhibits that the underlying field does not have much 
impact on the runtime or the result. The larger the dimension of
the considered Lie algebra is, the more superior is Affine. It
yields small dimensional representations and is the fastest of
all methods.

\subsection{Free nilpotent Lie algebras}

Next we consider the free nilpotent Lie algebras with $n$ generators of
class $c$ over the field $\F$, denoted $N_{n,c}(\F)$. 

\vspace*{0.2cm}
\begin{table}[htb]
\begin{center}
\begin{tabular}{|c|c|c|c|c|c|c|c|c|c|}
\hline
$n$ & $c$ & $\F$ & $\dim(N_{n,c}(\F))$ 
  & \multicolumn{2}{c|}{Regular} 
  & \multicolumn{2}{c|}{Dual} 
  & \multicolumn{2}{c|}{Affine}\\
\hline
& & & & time & $\dim$ & time & $\dim$ & time & $\dim$\\
\hline
2 & 5 & $\Q$ & 14 & 0.2 & 20 & 0.3 & 20 & 0.5 & 15\\
2 & 6 & $\Q$ & 23 & 0.9 & 34 & 1.3 & 34 & 8.4 & 24 \\
2 & 7 & $\Q$ & 41 & 3.2 & 65 & 4.8 & 65 & $\frownie$ &  $\frownie$ \\
2 & 8 & $\Q$ & 71 & 14 & 117 & 21 & 117 &  $\frownie$ & $\frownie$ \\
3 & 4 & $\Q$ & 32 & 0.8 & 41 & 1.7 & 41 & 54 & 33 \\
3 & 5 & $\Q$ & 80 & 11.5 & 113 & 17.5 & 113 &  $\frownie$ &  $\frownie$ \\
4 & 3 & $\Q$ & 30 & 0.9 & 36 & 1.3 & 36 & 37 & 31 \\
4 & 4 & $\Q$ & 90 & 13 & 113 & 19.7 & 113 &  $\frownie$ &  $\frownie$ \\
\hline
\end{tabular}
\vspace*{0.2cm}
\caption{Running times (in seconds) for $N_{n,c}(\Q)$.}\label{Table2}
\end{center}
\end{table}

Table $2$ displays the time in seconds for the three algorithms, with input 
$N_{n,c}$. The  $\frownie$ in the last two columns indicates that the 
algorithm Affine did not succeed, either because it made the ``wrong''
choice at some stage, or due to Memory problems: for its cohomology
computation it has to solve a system of linear equations which is of the
size $O(\dim(\Lg)^2)$ and this can be time-and space consuming.

\subsection{The Lie algebras $\Lf_n$}

Finally, we consider the Lie algebras $\Lf_n$ of the previous section. 
The results of that are contained in Table \ref{Table3}.

\vspace*{0.2cm}
\begin{table}[htb]
\begin{center}
\begin{tabular}{|c|c|c|c|c|c|c|c|}
\hline
$n$ & \multicolumn{2}{c|}{Regular} & \multicolumn{2}{c|}{Quotient} & 
\multicolumn{2}{c|}{Dual} & Affine \\
\hline
& time & $\dim$ & time & $\dim$ & time & $\dim$ & \\
\hline
13 & 8.6 & 85 & 14 & 43 & 12.3 & 43 & $\frownie$\\
14 & 17 & 105 & 28 & 53 & 24.7 & 53 &  $\frownie$\\
15 & 33 & 145 & 63 & 64 & 50 & 64  &  $\frownie$\\
16 & 64 & 185 & 125 & 77 & 102 & 77 &  $\frownie$\\
17 & 123 & 256 & 323 & 94 & 218 & 94 &  $\frownie$ \\
18 & 234 & 316 & 731 & 111 & 461 & 111 &  $\frownie$\\
19 & 487 & 433 & 1844 & 134 & 1162 & 134 &  $\frownie$\\
20 & 920 & 538 & 4009 & 158 & 3039 & 158&  $\frownie$\\
\hline
\end{tabular}
\vspace*{0.2cm}
\caption{Running time (in seconds) for the Lie algebras $\Lf_n$.}\label{Table3}
\end{center}
\end{table}

Table $3$ displays the time in seconds for the algorithms Quotient and 
Dual, with input $\Lf_n$. The $\frownie$ in the last column indicates 
that the algorithm Affine did not succeed. In this case, this was due
to the fact that Affine did not find any possible faithful representation
of dimension $n+1$. Of course, if our conjecture on $\Lf_n$ holds, then
it cannot succeed.

Note that the dimensionals of the determined modules for $\Lf_n$ are 
significantly larger than $n+1$. However, they do not seem to grow very 
fast. Some naive tests with least squares fits seem to suggest that 
the dimensions grow quadratically or cubically.

\subsection{Some comments}

From the above tables we conclude that if Affine succeeds, then it usually 
finds a module of significantly smaller dimension than Regular, Quotient or
Dual. This supports the suggested strategy to try this algorithm first.


\end{document}